\documentclass[12pt]{article}
\usepackage{amssymb}

\def\C{\mathbb C}

\def\N{\mathbb N}
\def\Z{\mathbb Z}

\newtheorem{thm}{Theorem}[section]
\newtheorem{lem}{Lemma}[section]
\newtheorem{defn}{Definitions}[section]

\begin{document}

\sffamily

\title{Derivatives of meromorphic functions of finite order}
\author{J.K. Langley}
\maketitle

\begin{abstract}
A result is proved concerning meromorphic functions $f$ of finite order in the plane such that all but finitely many
zeros of $f''$ are zeros
of $f'$. 
A.M.S. MSC 2000: 30D35.

\end{abstract}

\section{Introduction}\label{intro}

The starting point of this paper is the following theorem from \cite{Lalogsing}.

\begin{thm}[\cite{Lalogsing}]\label{thmA}
Assume that the function
$f$ is meromorphic of finite lower
order in the plane and that $f^{(k)} $ has finitely many zeros, for
some $k \geq 2$. Assume further that there exists a positive real number $M$
such that if $\zeta $ is a pole of $f$ of multiplicity $m_\zeta$ then
\begin{equation}
m_\zeta \leq M + |\zeta|^{M} .
\label{multest}
\end{equation}
Then $f$ has 
finitely many poles.
\end{thm}
Condition (\ref{multest}) is evidently satisfied if $f$ has finite order. 
Theorem \ref{thmA} fails for $k=1$, as shown by simple examples,
and for $k \geq 2$ and infinite lower order, in which case an example
is constructed in \cite{Lanew} 
with 
infinitely many poles, all simple, such that $f^{(k)}$ 
has no zeros at all. The result was inspired by the conjecture made by A.A. Gol'dberg, to the effect that for $k \geq 2$ and a
meromorphic function $f$ in the plane, regardless of growth, 
the
frequency of distinct poles of $f$ is controlled by the frequency of zeros of $f^{(k)}$, up to an error term which is small compared to the 
Nevanlinna characteristic.
Yamanoi has now proved this conjecture in a landmark  paper \cite{yamanoi}; however, because of the error terms involved,
his  result does not imply 
Theorem~\ref{thmA} directly.

This paper is concerned with a generalisation of Theorem~\ref{thmA} in a different direction.  The assumption there
that $f^{(k)} $  has finitely many zeros is a strong one, so that  it is natural to ask
whether it may be replaced by something less restrictive. A reasonable candidate is 
the condition
that all but finitely many zeros of $f^{(k)}$ have the same image under $f^{(k-1)}$, which may then be assumed to be~$0$, but
the following example shows that this does not  by itself imply that $f$ has finitely many poles.  
Set
\begin{equation}
f(z) = z- \tan z , \quad f'(z) = 1- \sec^2 z  = - \tan^2 z, \quad f''(z) = - 2 \tan z \sec^2 z .
\label{extan}
\end{equation}
Here all zeros of $f''$ are zeros of $f'$ and fixpoints of $f$, all zeros and poles of $f'$ have the same multiplicity,
and $1$ is an asymptotic value of $f'$. 
More generally it may be observed that, for any even positive integer $n$, the antiderivative
of $\tan^n z$ is meromorphic in $\C$. 
The example (\ref{extan}) shows that the following theorem, which evidently  implies Theorem \ref{thmA}, is essentially sharp.

\begin{thm}\label{thm1a}
Let $k \geq 2$ be an integer and let $f$ be a meromorphic function of finite lower
order in the plane with the following properties:
\\
(i) the zeros of $f^{(k-1)}$ have bounded multiplicities;\\
(ii) all but finitely many zeros  of $f^{(k)}$ are zeros of $f^{(k-1)}$;\\
(iii) there exists  $M \in (0, + \infty) $
such that if $\zeta $ is a pole of $f$ of multiplicity $m_\zeta$ then (\ref{multest}) holds;\\
(iv) for each $\varepsilon > 0$, all but finitely many zeros $z$ of $f^{(k)}$ satisfy
either $|f^{(k-2)}(z)| \leq \varepsilon |z|$ or $\varepsilon |f^{(k-2)}(z)| \geq |z|$.
\\
Then $f^{(k)} $ has a representation $f^{(k)} = Re^P$ with $R$ a rational function and $P$ a polynomial. In particular,
$f$ has finite order and finitely many poles, and $f^{(k)}$ has finitely many zeros.
\end{thm}
It  suffices to prove Theorem \ref{thm1a} for $k=2$ and, as already noted,
condition (iii) holds when $f$ has finite order. 
If $f$ is a meromorphic function of finite lower
order in the plane satisfying condition (ii)  of Theorem~\ref{thm1a}, with $k=2$,
then $f'$ has finitely many 
critical values and so finitely many asymptotic values,  by 
a result of Bergweiler and Eremenko \cite{BE} and its extension by Hinchliffe \cite{hinchliffe} to functions of finite lower order
(see Section \ref{critical2}).
Therefore Theorem \ref{thm1a} follows from the next result, which fails for infinite lower order,
because of the same example from \cite{Lanew} mentioned after Theorem \ref{thmA}.

\begin{thm}\label{thm1}
Let $f$ be a meromorphic function of finite lower
order in the plane satisfying conditions (i), (ii) and (iii) of Theorem \ref{thm1a}, with $k=2$. 
Assume  that there exist positive real numbers $\kappa$ and $R_0$ such that if $z$ is a zero of $f''$ 
with $|z| \geq R_0$ then  
$|f(z) - \alpha z | \geq \kappa |z|$ for all  finite non-zero
asymptotic values $\alpha$ of $f'$. 
Then $f'' = Re^P$ with $R$ a rational function and $P$ a polynomial.
\end{thm}

\section{Lemmas needed for Theorem \ref{thm1}}\label{lemmas}

Throughout this paper
$B(z_0, r)$ will denote the disc $\{ z \in \C : | z - z_0 | < r \}$ and
$S(z_0, r)$ will be the circle $\{ z \in \C  : | z - z_0 | = r \}$.
The following results are both well known. 

\begin{lem}[\cite{Tsuji}, p.116]\label{lemC}
Let $D$ be a
simply connected domain not containing the origin, and let
$z_0$ lie in $D$. Let $r$ satisfy $0 < 4 r < | z_0| $ or
$4 | z_0| < r < \infty $.
Let $\theta(t)$ denote 
the angular measure of $D \cap S(0, t)$,
and
let $D_r$ be the component of $D \setminus S(0, r)$
which contains $z_0$.
Then the harmonic measure of $S(0, r) $
with respect to the domain $D_r$, evaluated at $z_0$, satisfies
\begin{equation}
\omega (z_0, S(0, r) , D_r ) \leq C
\exp \left( - \pi   \int_I \frac{dt}{ t \theta(t) }  \right),
\label{tsujiest1}
\end{equation}
with $C$ an absolute constant, 
$I = [ 2 |z_0| , r/2 ]$ if $ 4 |z_0| < r$, and
$I = [ 2 r , |z_0|/2 ]$ if $4 r <  |z_0|$. 
\end{lem}


\begin{lem}[\cite{Hay7}, p.366]\label{cartan}
Let $Q$ be a positive integer and let
$w_1, \ldots , w_Q$ be complex numbers. For each $\Lambda > 0$ the estimate
\begin{equation}
 \prod_{j=1}^Q |z - w_j| \geq \Lambda^Q
\label{cartanest}
\end{equation}
holds for all $z$ 
outside a union of discs having sum of radii at most $6 \Lambda$.
\end{lem}


\section{Critical points and
asymptotic values }\label{critical2}

Suppose that the function $h$ is transcendental and meromorphic in the
plane, and that $h(z)$ tends to 
$a \in \C$ as $z$ tends to 
infinity along a path $\gamma$. Then $a$ is an asymptotic value of 
$h$, and the inverse
function $h^{-1}$ has a
transcendental singularity over $a$ \cite{BE, Nev}.
For each $t > 0$, let
$C(t)$ be
that component of $C'(t) = \{ z \in \C : | h(z) - a | < t \}$ which
contains an unbounded subpath of~$\gamma$. 
The singularity of $h^{-1}$ over $a$ corresponding to $\gamma$
is called direct \cite{BE} if $C(t)$, for some $t > 0$,
contains no zeros
of $h(z) - a$.
Singularities over
$\infty$ are classified 
analogously.

Recall next
some standard facts from \cite[p.287]{Nev}.
Suppose that $G$ is a transcendental meromorphic function
with
no asymptotic or
critical values in $1 < |w| < \infty $.
Then every
component $C_0$ of the set
$\{ z \in \C : |G(z)| > 1 \}$ is simply connected, and there are
two possibilities. Either (i) $C_0$ contains 
one pole $z_0$ of $G$ of multiplicity $k$, in which case 
$G^{-1/k} $ maps $C_0$ univalently onto $B(0, 1)$, or 
(ii) $C_0$ contains no pole of $G$, but instead a path
tending to infinity on which $G$ tends to infinity. In case (ii) the function $w = \log G(z)$ maps 
$C_0$ univalently onto the right half plane.

\begin{lem}[\cite{Lalogsing}]\label{lem1}
There exists a positive absolute constant $C$ with the following property. 
Suppose that $G$ is a transcendental meromorphic function in the plane
and that $ G' $ has
no asymptotic or critical values $w$ with
$0 < |w| < d_1 < \infty $.
Let
$D$ be a component of the set
$\{ z \in \C : |G'(z)| < d_1 \}$ on which $G'$ has no zeros, but such that
$D$
contains a path tending to infinity
on which $G'(z)$ tends to $ 0$.
If $z_1$ is in $D$ and $\log |d_1/G'(z_1)| \geq 1$ then
$$
|G(z_1)| \leq S + \frac{C |z_1G'(z_1)|}{ \log |d_1/G'(z_1)| }  ,
$$
in which the
positive constant $S$ depends on $G$ and $D$ but not on $z_1$. 
\end{lem}
\hfill$\Box$
\vspace{.1in}

Suppose next that the function $F$ is meromorphic of finite lower order in the plane,
and that all but finitely many zeros of 
$F'$ are zeros of $F$. Then $F$ has finitely many critical values.
By Hinchliffe's extension \cite{hinchliffe} to the finite lower order case
of a theorem of Bergweiler and Eremenko \cite{BE}, the function $F$ 
has finitely many asymptotic values.
Furthermore, 
all asymptotic values of $F$
give rise to direct
transcendental singularities of the inverse function $F^{-1}$ 
and,
by the Denjoy-Carleman-Ahlfors theorem \cite {BE, Hay7, Nev},
there are finitely many such singularities.
The following facts are related to the argument from \cite[Section 4]{Lanew}.
Let $J$ be a  polygonal Jordan curve in $\C \setminus \{ 0 \}$ such that every
finite non-zero critical or asymptotic value of $F$ 
lies on $J$, but is not a vertex of $J$, and such that 
the complement of $J$ in $\C \cup \{ \infty \} $
consists of two simply connected domains $B_1$ and $B_2$, with
$0 \in B_1$ and $\infty \in B_2$. Fix conformal
mappings 
\begin{equation}
h_m : B_m \to  B(0, 1), 
\quad m = 1, 2, \quad h_1(0)=0, \quad h_2( \infty ) = 0.
\label{a1}
\end{equation}
The mapping $h_1$ may then be extended to be quasiconformal on the plane, fixing infinity, and there
exist a meromorphic function $G$  and a quasiconformal mapping $\psi$ such that 
$h_1 \circ F = G \circ \psi $ on $\C$. It follows that for $j=1,2$ all components of $F^{-1}(B_j)$ are simply connected  and all but finitely
many  are unbounded, since
all but finitely many zeros $z$ of $G'$ have $G(z) = 0$. 

\section{Proof of Theorem \ref{thm1}: first part}\label{proof1}

Let the function $f$ be as in the hypotheses. If $f''/f'$ is a rational function then $f'$ is a rational function multiplied by the 
exponential of a polynomial,
and so is $f''$. Assume henceforth that $f''/f'$ is transcendental: then obviously so is $f$.  
Apply the reasoning and notation of Section~\ref{critical2}, with 
$F = f'$.
The following  is an immediate
consequence of Lemma \ref{lem1}.

\begin{lem}\label{lem3}
Arbitrarily small positive real numbers
$\varepsilon_1$ and $\varepsilon_2$ may be chosen with the following properties.
There exist finitely many 
unbounded simply connected domains $U_n$, each of which is a component of the set
$\{ z \in \C : | f'(z) - b_n | < \varepsilon_1 \}$,
such that $U_n$ 
contains a
path tending to infinity on which $f'(z)$ tends to the  finite asymptotic value $b_n $. Here
$f'(z) \not = b_n$ on $U_n$ and
$|f(z) - b_n z| < \varepsilon_2 |z|$ for all
large $z$ in $U_n$.
If $\Gamma$ is a path tending to infinity on which
$f'$ tends to a finite asymptotic value $\alpha$, then there exists $n$ such that
$\alpha = b_n$ and $\Gamma \setminus U_n$ is bounded. 
The  $b_n$ need not be distinct, and some of them may be~$0$.
\end{lem}
\hfill$\Box$
\vspace{.1in}
\begin{lem}
 \label{lemexth1}
There exists a positive real number  $s_1 < \varepsilon_1$ with the following property. 
Let $b_p$ be a finite non-zero asymptotic value of $f'$. Then 
the conformal map $h_1: B_1 \to B(0, 1)$ extends to be analytic and univalent on
$B_1 \cup B(b_p, s_1)$.
\end{lem}
\textit{Proof.} This follows from the Schwarz reflection principle and the fact that each non-zero
$b_p$ lies on the polygonal Jordan curve $J = \partial B_1 $ but is not a vertex of $J$. 
\hfill$\Box$
\vspace{.1in}

\begin{defn}
 \label{def1}
Fix  positive real numbers $\rho$, $\sigma $ and $ \tau$
with $  \tau < s_1 < \varepsilon_1$
and $\sigma/ \tau$ and 
$\rho/\sigma$ small. 
Fix $W_0 \in \C$ such that $f'(W_0) $ is large.
\end{defn}

\begin{lem}
 \label{newlem0}
With the notation of Definitions \ref{def1}, there exist positive real numbers $M_1$, $M_2$, $M_3$ having the following properties.
Let $z_0$ be large 
with $ |f'(z_0)| < \tau $ and assume that $z_0$ lies in a component
$C$ of $(f')^{-1}(B_1)$ satisfying one of the following two conditions:\\
(A) there is  at least one zero of $f''$ in $C$; \\
(B) the function $f'$ is univalent on $C$, and $C \cap U_p$ and $C \cap U_q$ are both non-empty, where $U_p$ and $U_q$ are 
 as in Lemma \ref{lem3} with  $0 \neq b_p \neq b_q \neq 0$.\\
Then $|z_0 f''(z_0) | \leq M_1$ and there exists a disc 
$B( z_0^*, M_2 |  z_0^* | ) \subseteq B\left(z_0, \frac12 |z_0|\right) \cap  C$ on which
\begin{equation}
\left| \frac{f'''(\zeta)}{f''(\zeta)} \right| \leq \frac{M_3}{|\zeta|}  .
\label{lem10a}
\end{equation}
\end{lem}
\textit{Proof.} 
Observe that conditions (A) and (B) are mutually exclusive.
Denote positive constants by $c_j$ and small positive constants by $\delta_j$; these will
be independent of $z_0$ and $C$. In case (A) there is exactly
one point  in $C$ at which $f''$ vanishes, and it must be a zero of $f'$. In both cases $f'(C) = B_1$
(see Section~\ref{critical2}), and 
$C$ contains precisely one zero $z_1$ of $f'$, of multiplicity $ m \leq c_1$, by hypothesis (i)
of the theorem, with $m=1$ in Case (B).
There exist only finitely many components $C_1$ of $(f')^{-1}(B_1)$ which are bounded or have a zero of $f''$ 
on their boundary,
and if one of these contains a zero of $f'$ 
then the set $\{ z \in C_1 : |f'(z)| \leq \tau \}$ is compact.
Therefore since $z_0$ is large the component $C$ is unbounded and simply connected and its boundary $\partial C$ contains no zeros of
$f''$. Now set
$v_0 = (h_1 \circ f')^{1/m} $, with $h_1$ as in (\ref{a1}). Then $v_0$
maps $C$ univalently onto  $B(0, 1)$, and $u_0 = v_0( z_0)$ satisfies
$|u_0| \leq \delta_1$, since $m \leq c_1$ and $\tau$ is small. 

Let $\Gamma$ be a component of $\partial C$. Then $\Gamma$ is a simple curve tending to infinity in both directions
and, as $z$ tends to infinity in either direction along $\Gamma$, the image $f'(z)$ must tend to a finite non-zero
asymptotic value of $f'$; this is because $v_0$ is univalent on $C$. 
Hence there exists $z_1$ lying close to $\Gamma$, such that
$z_1 \in C \cap U_n$, for some $U_n$ as in Lemma \ref{lem3}, with $b_n \neq 0$
and $|f'(z_1) - b_n| < \varepsilon_1 $. 
By construction, $b_n$ lies on the polygonal Jordan curve $J$ but is not a vertex of $J$. Thus
analytic continuation of $(f')^{-1}$ along a path in the semi-disc $B(b_n, \varepsilon_1 ) \cap B_1$ then gives a point
$z_2 \in C \cap U_n$ with 
$|f'(z_2) - b_n| < \varepsilon_1 $, as well as 
$|h_1 (f'(z_2))| \leq 1 - \delta_2$, which implies in turn that $|v_0(z_2)| \leq 1 - \delta_3$. 

Let $G_0 : B(0, 1) \to C$ be the inverse function of $v_0$, and suppose that 
$G_0'(u_0) = o( |z_0| )$. Then Koebe's distortion theorem implies that
$G_0'(u) =  o( |z_0| )$ for $|u| \leq 1 - \delta_3$. In Case (A) this gives a path $\gamma$
in $C$, of length $o( |z_0| )$, joining $z_3 = G_0(0)$ to $z_2$ via $z_0$, and with $|f'(z)| \leq c_2$ on~$\gamma$.
Since $z_0$ is large
so are $z_2$ and $z_3$. Thus Lemma \ref{lem3} and integration of $f'$ yield
\begin{equation}
f(z_3) = f(z_2) + o( |z_0|), \quad |f(z_3) - b_n z_3 | \leq \varepsilon_2 |z_2|  +  o( |z_0| )
\leq ( \varepsilon_2 + o(1)) |z_3| .
\label{z3est}
\end{equation}
But by the assumption of Case (A),
$f''$ has a zero in $C$, which must be at $z_3$, so that,  by the hypotheses
of the theorem,  $|f(z_3) - b_n z_3 | \geq \kappa | z_3 |$. This contradicts (\ref{z3est}),
since $ \varepsilon_2 $ is small.
Next, in Case (B) the above analysis may be applied twice, to give a path $\gamma$ in $C$ of length $o( |z_0| )$, on which
$|f'(z)| \leq c_2$,  such that $\gamma$ joins points $w_p \in C \cap U_p$ and $w_q \in C \cap U_q$ via $z_0$, where 
$b_p$ and $b_q$ are distinct and non-zero, and $|f'(w_j) - b_j| <  \varepsilon_1$ for $j = p, q$. Therefore the $w_j$ satisfy
$w_j \sim z_0$ and 
$|f(w_j) - b_j w_j| \leq  \varepsilon_2 |w_j| \leq 2  \varepsilon_2 | z_0|$ for $j = p, q$. 
Since $ \varepsilon_2$ is small and 
integration of $f'$ along $\gamma$ leads to $f(w_p) - f(w_q) = o( |z_0| )$, this case also delivers a contradiction.

It follows in both cases that $|G_0'(u_0)| \geq c_3 |z_0|$, which implies at once that  $| z_0 v_0'(z_0) | \leq c_4$. Writing
$f'(z) = h_1^{-1} (v_0(z)^m)$ and using the fact that $ |f'(z_0)| < \tau $ and $m \leq c_1$ gives  $| z_0 f''(z_0) | \leq c_5$.
To prove the last assertion requires a disc on which $f'$ is univalent.
To this end, observe that $|G_0'(u_0)| \leq c_6 |z_0|$, since $z_0$ is large but $C$ does not contain the point
$W_0$ chosen in Definitions \ref{def1}. 
Now choose $u_0^*$ with $|u_0^* - u_0| \leq \delta_4 $ and $|u_0^*| \geq \delta_4$, and 
choose $\delta_5$ so small that the function $u^m$ is univalent on $B(  u_0^* , \delta_5  )$.
Then
Koebe's distortion theorem 
implies that the image $X_0$ of
$B(  u_0^* , \delta_5  )$ under $G_0$ lies in $ B\left(z_0, \frac12 |z_0|\right) \cap  C$
and 
contains a disc $B( z_0^* , 2 M_2 |  z_0^* | )$, where $ z_0^* = G_0 ( u_0^* )$ and
$M_2 = \delta_6$: this requires only that $\delta_4$ and $\delta_5/\delta_4$ be small enough, independent of $z_0$.
The function $v_0(z)^m$ is univalent on $X_0$ and therefore so
is $f'$. 
Now take $\zeta$ in $B( z_0^*,  M_2 |  z_0^* | )$ and set
$$
g(z) = \frac{ f'(\zeta+ M_2 |  z_0^* | z) - f'(\zeta)}{M_2 |  z_0^* | f''(\zeta) } = z + \sum_{\mu=2}^\infty A_\mu z^\mu 
$$
for $|z| < 1$, so that the estimate (\ref{lem10a}) follows from Bieberbach's bound $|A_2| \leq 2$. 
\hfill$\Box$
\vspace{.1in}

It will be seen that hypothesis (i) of Theorem \ref{thm1} plays a key role in the above proof of Lemma \ref{newlem0},
principally by preventing $z_0$ from lying too close to the boundary of $C$. 



\begin{lem}\label{newlem1}
With the notation of Lemma \ref{lem3} and Definitions \ref{def1}, let $z_1$ be large and satisfy
\begin{equation}
z_1 \in U_p, \quad b_p \neq 0, \quad  \sigma < |f'(z_1) - b_p| < \tau  < s_1 , \quad f'(z_1) \in J = \partial B_1 ,
\label{z1def}
\end{equation}
and let $C$ be the component of $(f')^{-1}(B_1)$ with $z_1 \in \partial C$.
Assume that one of the following two mutually exclusive conditions
holds:\\
(a) the function $f'$ is not univalent on $C$;\\
(b) the function $f'$ is univalent on $C$, and $C \cap U_q$ is non-empty, for some $q$ with $0 \neq b_q \neq b_p$. \\
Then there exists an open set $H_1$, with
\begin{equation}
 \label{z1*def}
H_1 \subseteq B\left(z_1, \frac12 |z_1|\right) \cap C 
\quad \hbox{and} \quad \partial H_1 \cap \partial C = \{ z_1 \} ,
\end{equation}
such that $f'$ maps $H_1$ onto an open disc $K_1 \subseteq B_1$, of diameter less than $\rho$,
which is tangent to $J = \partial B_1 $ at
$ f'(z_1)$. Furthermore, 
$H_1$ contains an open disc $L_1$ of radius $M_4 |z_1| $ on which (\ref{lem10a}) holds; here both $M_3$ and $M_4$ are independent of $z_1$
and $C$. 
\end{lem}
\textit{Proof.} 
The component $C$ is unique because $z_1$ is large and $f''$ has finitely many zeros which are not zeros of $f'$. 
As in Lemma \ref{newlem0} denote small positive constants by $\delta_j$, and positive constants by $c_j$; these will again be
independent of $z_1$ and $C$. 
Let $\gamma_0 $ be the straight line segment 
$$
u = tu_1, \quad \delta_1 \leq t \leq 1 , \quad u_1 = h_1(f'(z_1)) \in S(0, 1),
$$
where $\delta_1$ is chosen sufficiently small that
$|h_1(w)| \leq \delta_1 $ implies that $|w| \leq \delta_2 < \tau <  \varepsilon_1$.
Using (\ref{z1def}) and
the conformal extension of $h_1$ to $B_1 \cup B(b_p, s_1)$ given by Lemma \ref{lemexth1}, 
define 
domains $F_1 \subseteq B_1 \cup \{ \zeta \in \C: \rho < | \zeta - b_p | < s_1 \}$
and $E_1$ by 
$$
E_1 = \{ u \in \C : {\rm dist} \, \{ u, \gamma_0 \} < \delta_3 \} = h_1 (F_1),
$$
in which $\delta_3$ is small compared to $\delta_1$, which ensures that $0 \not \in E_1$. 
Then $F_1$ contains no singular values of the inverse function $(f')^{-1}$, and $z_1$ lies in a component $D$ of 
$(f')^{-1}(F_1)$ such that $h_1 \circ f'$ maps $D$ conformally onto $E_1$.
Let $G_1: E_1 \to D$ be the inverse function of $ h_1 \circ f'$, 
and choose $z_2 \in D$ with
$u_2 = h_1(f'(z_2)) = \delta_1 u_1$ and hence $  |f'(z_2)| \leq  \delta_2 < \tau$. 
Observe that $z_2$ lies in $C$. 
Repeated application of the Koebe distortion theorem yields 
$
c_1 | G_1'(u_1)| \leq |G_1'(u)| \leq c_2 | G_1'(u_1)| 
$
on
the line segment 
$ \gamma_0$, 
and the image $\sigma_1 = G_1( \gamma_0) $ is a path of length at most $c_3 | G_1'(u_1)|   $ from $z_1$ to $z_2$ in~$D$. 

Suppose first that $G_1'(u_1) = o( |z_1|)$. Then  $z_2 \sim z_1 $ and $G_1'(u_2) =  o( |z_1|)$,
from which it follows that
$z_2 f''(z_2)$ is large.  Hence  $C$ satisfies neither 
condition (A) nor condition (B) of Lemma \ref{newlem0}, and so cannot satisfy (b), because (b) implies (B)
since $C \cap U_p \neq \emptyset$ and $b_p \neq 0$. 
Hence $f'$ is not univalent on $C$ but $C$ contains no zero of $f''$.
Thus $C$ must contain a path $\Gamma$ tending to infinity on which $f'(z)$ tends to $0$, and $C$ meets 
one of the components 
$U_n$ with $b_n = 0$. Moreover, $\log (h_1 \circ f')$ maps $C$ univalently
onto the left half plane (see Section~\ref{critical2}). Therefore, since $  |h_1(f'(z_2))| \leq  \delta_1 $, there exists a path 
$\Gamma'$ in $C$ joining $z_2$ to some
$z_3 \in \Gamma$ on which $|h_1(f'(z))| \leq \delta_1$ and $|f'(z)| < \varepsilon_1$, and hence $z_2 \in U_n$. 
Since $z_1$ is large, and $z_2 \sim z_1 $, 
Lemma \ref{lem3} gives 
$|f(z_2) | \leq \varepsilon_2 |z_2| $ and $|f(z_1) - b_p z_1 | \leq \varepsilon_2 |z_1|$, in which $b_p \neq 0$.
On the other hand $|f'(z)| \leq c_4$ on $\sigma_1$, and so integration yields  
$f(z_1) = f(z_2) + o( |z_1|) $
and a contradiction.

It must therefore be the case that $| G_1'(u_1)| \geq  c_5 |z_1|$. 
However, the point
$W_0 $ chosen in Definitions \ref{def1} is not in $D$ and so $| G_1'(u_1)| \leq  c_6 |z_1|$.
Now let $G_2 = G_1 \circ h_1 : F_1 \to D$ be the inverse function of $f'$, and set $v_1 = f'(z_1) = h_1^{-1}(u_1) \in J$. 
Then (\ref{z1def}) yields 
$c_7 |z_1| \leq |G_2'(v_1)| \leq c_8 |z_1|$, as well as 
$B(v_1, 2 \delta_4 ) \subseteq F_1$ for some $\delta_4 < \rho $, and  
Koebe's distortion theorem gives
$ c_9 |z_1| \leq | G_2'(v)| \leq  c_{10} |z_1|$ on $B(v_1,  \delta_4 ) $. Hence
$G_2(  B(v_1 ,  \delta_5 ) ) \subseteq B\left(z_1, \frac12 |z_1| \right)$, provided
$\delta_5 \leq \delta_4$ is chosen small enough. Let $K_1 \subseteq B(v_1 ,  \delta_5 ) \cap B_1$ be an open disc of radius
$\delta_6 \leq \frac14 \delta_5 $, which is tangent to $J$ at $v_1$. Then $H_1 = G_2(K_1) $ satisfies (\ref{z1*def}), and $H_1$ contains
a disc $B(z_1^*, 2 M_4   |z_1| )$, with $M_4 = \delta_7$. It may now be assumed that $M_3$ is large enough that
(\ref{lem10a}) 
holds on $L_1 = B(z_1^*,  M_4   |z_1| )$, since  Bieberbach's theorem may be applied as in the proof of Lemma \ref{newlem0}.
\hfill$\Box$
\vspace{.1in}

\section{The frequency of poles 
of $f$ 
and zeros of $f''$}\label{frequency}

\begin{lem}
 \label{poleslem}
Let $w_1, \ldots , w_Q$ be pairwise distinct  poles of $f$ with $|w_j|$ large.
For  $1 \leq j \leq Q$ let $D_j$ be
the component of $(f')^{-1}(B_2)$ in which $w_j$ lies.
Then for each $j$ there exists $p_j \in \Z$ such that $\partial D_j$ contains a Jordan arc $\lambda_j$ which is mapped univalently by $f'$ onto 
a line segment $\mu_j$ of length at least $\sigma$, and these may be chosen so that
\begin{equation}
 \label{mujdefn}
 \lambda_j \subseteq U_{p_j}, \quad 
\mu_j \subseteq  
\{ \zeta \in J = \partial B_2 : \sigma < | \zeta - b_{p_j} | < \tau \} , \quad  b_{p_j} \neq 0, 
\end{equation}
where $U_{p_j}$ and $b_{p_j}$ are as in Lemma \ref{lem3}, while $\sigma$ and $\tau$  are as in  Definitions \ref{def1}. 

Moreover, if points $z_j$ are chosen such that $z_j \in \lambda_j$ for  $1 \leq j \leq Q$, then each $|z_j|$ is large and
for each $j$ there exists an open disc $L_j \subseteq  B\left(z_j, \frac12 |z_j| \right)$ of radius $M_4 |z_j| $, on which (\ref{lem10a}) holds,
where  $M_4$ is as in Lemma \ref{newlem1}.
The $L_j$ are pairwise disjoint. 
\end{lem}
\textit{Proof.} By the discussion in Section \ref{critical2}, each $D_j$ is unbounded and simply
connected and the boundary $\partial D_j$ contains no zeros of $f''$. 
Each component of $\partial D_j$ is a simple path tending to infinity in both directions, 
and there exists a component $\Gamma_j$ of $\partial D_j$ which separates $w_j$ from the point $W_0$ chosen 
in Definitions \ref{def1}. Since $D_j$ contains a pole of $f$ it follows that $f'$ is finite-valent
on $D_j$. 
Thus as $z$ tends to infinity in either direction along $\Gamma_j$ the image $f'(z)$ must tend to a non-zero finite asymptotic value  of
$f'$. 
In particular, $\Gamma_j$ meets some $U_p$ as in Lemma \ref{lem3} with $b_p \neq 0$, and following $\Gamma_j$ while staying in 
$U_p$ gives $\lambda_j$ and $\mu_j$ as in (\ref{mujdefn}).
Furthermore, each
$w_j$ is large and, for any $M_5 > 0$, the disc $B(0, M_5)$ meets only finitely many components of $(f')^{-1}(B_2)$, 
each of which contains at most one pole of $f$. Hence if $z_j \in \lambda_j$ then $z_j$ is large.

To prove the existence of the $L_j$, 
choose for each $j$ a component $E_j$ of $(f')^{-1}(B_1)$ with $\Gamma_j \subseteq \partial E_j$. Since
$\Gamma_j$ separates the pole $w_j$  of $f$ from $W_0$ it follows that $\Gamma_j$ is not the whole boundary $\partial E_j$. In particular, if
$f'$ is univalent
on $E_j$ then 
$\Gamma_j$ must meet components $U_p$ and $U_q$ with $ b_p $ and $b_q$ distinct and non-zero. 
Thus each of these components $E_j$ of $(f')^{-1}(B_1)$ satisfies one of the conditions (a), (b)
of Lemma \ref{newlem1}, which may now be applied with $z_1$ replaced by each~$z_j$. This gives open sets $H_j \subseteq 
B\left(z_j, \frac12 |z_j| \right) \cap E_j$, 
each containing an open disc 
$L_j$ of radius $M_4 |z_j| $ on which (\ref{lem10a}) holds. Moreover, $f'$ maps $H_j$ onto a disc $K_j \subseteq B_1$ which is tangent to $J$ at
$f'(z_j) $ and has diameter less than~$\rho$.

To show that the $L_j$ are disjoint, suppose that $1 \leq j < j' \leq Q$ and that $H_j \cap H_{j'} \neq \emptyset$, from which it follows of course
that $K_j \cap K_{j'} \neq \emptyset$. Since  $\rho$ is  small compared to $\sigma$ and $z_j \in \lambda_j$, the open disc $U = B(f'(z_j), 3 \rho)$ 
contains no singular value of $(f')^{-1}$, by (\ref{mujdefn}). 
But $K_j$ and $K_{j'}$ have diameter less than $\rho$, and so their closures lie in $U$. 
Thus $H_j$ and $H_{j'}$ both lie in the same component of $(f')^{-1}(U)$, as do $z_j$ and $z_{j'}$,
which forces  $\Gamma_j = \Gamma_{j'} $ and gives a contradiction. 
\hfill$\Box$
\vspace{.1in}

\begin{lem}\label{lem2}
Let
$L(r) \to \infty $ with $L(r) \leq \frac18 \log r$
as $r \to \infty $, and for $k > 0$ and large $r$ define the annulus $A(k)$ by
$A(k) = \{ z \in \C : r e^{-k L(r) } \leq |z| \leq re ^{k L(r) }  \}$.
Then the number
$N_1$ of distinct poles of $f$ and zeros of $f''$ in $A(1)$ satisfies
\begin{equation}
N_1 = O ( \phi (r) )  \quad \hbox{as $ r \to \infty $, where} \quad
\phi (r) =  L(r) + \frac{ \log r}{ L(r)} .
\label{lem2a}
\end{equation}
\end{lem}
\textit{Proof.} 
Assume that $r$ is large and that $A(1)$ contains 
$Q = 2N$ distinct poles $w_1, \ldots , w_{2N}$ of $f$,
with 
$\phi (r) = o( N)$. For $j=1, \ldots, Q$ let $D_j$ be
the component of $(f')^{-1}(B_2)$ in which $w_j$ lies, let $q_j$ be the multiplicity of the pole of $f'$ at
$w_j$.
Each $D_j$ is unbounded and simply connected and may be assumed not to contain the origin. 
Let $v_j = ( h_2 \circ f' )^{1/q_j}$, so that 
$v_j$ maps $D_j$ conformally onto $B(0, 1)$,
with $v_j(w_j) = 0$.

For $0 < t < \infty$ let $\theta_j(t)$ be the angular measure of
$D_j \cap S(0, t)$. 
Let $c$ denote positive constants, not
necessarily the same at each occurrence, but not depending
on $r, L(r)$ or $N$.
For $m \in \N$ the Cauchy-Schwarz inequality gives
$m^2 \leq 2 \pi \sum_{j=1}^m 1/ \theta_j (t) $ so that, as in
\cite{Lanew},
at least $N$ of the $D_j$ have
\begin{equation}
\int_{2r e^{L(r)} }^{(1/2)r e^{2L(r)} } 
\, \frac{dt}{ t \theta_j (t) } > cN  L(r) ,
\quad
\int_{2r e^{- 2 L(r) }}^{(1/2)r e^{ - L(r) }} \, \frac{dt}{ t \theta_j (t) } > cN  L(r) .
\label{b2}
\end{equation}
It may be assumed after re-labelling if necessary that (\ref{b2}) holds for $D_1 , \ldots , D_{N} $. 
Since $w_j$ lies in
$A(1)$, it follows from Lemma \ref{lemC} that
$$
\omega(w_j, \sigma_j, D_j) 
\leq
c \exp \left( -  \pi \int_{2r e^{L(r)}}^{(1/2)r
e^{2L(r)}} \, \frac{dt}{ t \theta_j (t) }\right) + 
c \exp \left( -   \pi 
\int_{2r e^{- 2L(r)}}^{(1/2)r e^{ - L(r) }} \, \frac{dt}{ t \theta_j (t) } \right)  .
$$
Combining this  with (\ref{lem2a}),  (\ref{b2}) and condition (iii) of the theorem shows that
$\omega(w_j, \sigma_j, D_j) = o( 1/q_j)$
for $j = 1, \ldots , N$, where
$\sigma_j =   \partial D_j \setminus A(2)$. 
But Lemma \ref{poleslem} gives an arc $\lambda_j \subseteq \partial D_j$,
mapped by $f'$ onto a line segment $\mu_j \subseteq J$ as in (\ref{mujdefn}), of length at least $\sigma$.
Since $b_{p_j}$ in (\ref{mujdefn}) is not a vertex of $J$, while $\tau$ is small,
an application of the Schwarz reflection principle to $h_2$ shows  that $h_2 \circ f'$
maps $\lambda_j$ to an arc of $S(0, 1)$ of length at least $c$, and $v_j( \lambda_j)$ 
has
angular measure at least $c/q_j$. 
The conformal invariance of harmonic measure under $v_j$  implies that
$\lambda_j $ cannot be contained in $\sigma_j$, and so there exists $z_j \in \lambda_j \cap A(2)$. 
The corresponding $N$ pairwise disjoint discs $L_j$ given by Lemma \ref{poleslem} lie in 
the annulus $A(3)$, and hence
$$
cN \leq \sum_{j=1}^N  \int_{L_j} |z|^{-2} dx dy \leq
\int_{A(3) } |z|^{-2} dx dy \leq c L(r) \leq c \phi (r) = o(N).
$$
This is a contradiction and the asserted upper bound for the number of distinct poles in $A(1)$ is proved. 
The same upper bound for the number of distinct zeros $\zeta_j$ of $f''$ in $A(1)$ 
follows at once from Lemma \ref{newlem0},
because such zeros give rise to pairwise disjoint discs $B( \zeta_j^*, M_2 | 
\zeta_j^*|) \subseteq A(2)$.
\hfill$\Box$
\vspace{.1in}

Since all but finitely many zeros of $f''$ are zeros of $f'$, which have bounded multiplicities by assumption,
choosing $L(r) = \frac18 \log r$ in Lemma \ref{lem2} gives
$$
\overline{n} ( r^{9/8} , f) -
\overline{n} ( r^{7/8} , f) + n( r^{9/8} , 1/f'') -
n ( r^{7/8} , 1/f'') = O( \log r ),
$$
and so
\begin{equation}
\overline{N}(r, f) + N(r, 1/f'') = O( \log r )^2 \quad \hbox{as}  \quad r \to \infty .
\label{polesfreq}
\end{equation}

\begin{lem}\label{lem00}
The lower order of $f''/f'$ is at least $\frac12$.
\end{lem}
\textit{Proof.} If this is not the case then
the function
$f'/f''$ has finitely many poles and
is transcendental of lower order less than $\frac12$.
The $\cos \pi \lambda $ theorem 
\cite{BarK} now gives
$r_j \to + \infty $ such that
$ f''(z)/f'(z) = O( r_j^{-2} )$ on $ S(0, r_j ) $.
Moreover, the main result of 
\cite{LRW} 
gives a path $\gamma$ tending to infinity with
$$
\int_\gamma \left| \frac{f''(z)}{f'(z)} \right| \, |dz| < \infty .
$$
This implies that, as $z$ tends to infinity 
in the union of $\gamma$ and the $S(0, r_j)$, the image  $f'(z)$ tends to some $b_n$ as in Lemma \ref{lem3},
contradicting the fact that the
$U_n$ are simply connected.
\hfill$\Box$
\vspace{.1in}

\begin{lem}\label{lem30}
The function $f''$ has the form $f'' = \Pi_1/\Pi_2$, where $\Pi_1$ and $\Pi_2$ are entire 
such that $\Pi_2$ has finite order and $\Pi_1 \not \equiv 0$ has order $0$. 
Moreover, the lower order of $\Pi_2$ is at least $1/2$.
\end{lem}
\textit{Proof.} Using (\ref{multest}) and (\ref{polesfreq}) 
shows that $N(r, f'')$ has finite order and 
$N(r, 1/f'')$ has order $0$. Since $f''$ has finite lower order, this gives the asserted representation for 
$f''$.
On the other hand, Lemma \ref{lem00} implies that $f'$ has lower order at least
$1/2$ and so has $f''$, and hence so has $\Pi_2$. 
\hfill$\Box$
\vspace{.1in}

\begin{lem}\label{wvlem}
Let $h(z) = z f'''(z)/f''(z)$. 
For all
$s \geq 1$ lying outside a set $E_0$ of finite logarithmic measure, 
there exists $\zeta_s$ with
$| \zeta_s| = s$ and $ | h( \zeta_s) | > s^{1/3} $.
\end{lem}
\textit{Proof.} Take $\Pi_1$ and $\Pi_2$ as in Lemma \ref{lem30}. Applying the Wiman-Valiron theory \cite[Theorem 12]{Hay5} and 
standard estimates for logarithmic derivatives \cite{Gun2} makes it possible to write,
for $| \zeta_s| = s$ with
$|\Pi_2( \zeta_s)| = M(s, \Pi_2)$ and $s$ outside a set of finite logarithmic measure,
$$
\frac{f'''}{f''} = \frac{\Pi_1'}{\Pi_1} - \frac{\Pi_2'}{\Pi_2}, \quad
\left| \frac{\Pi_2'(\zeta_s)}{\Pi_2(\zeta_s)} \right| \sim \frac{\nu(s)}{s} , \quad 
\left| \frac{\Pi_1'(\zeta_s)}{\Pi_1(\zeta_s)} \right| \leq s^{-3/4}.
$$
Here $\nu(s)$ is the central index of $\Pi_2$ and has lower order at least $1/2$. 
\hfill$\Box$
\vspace{.1in}

\section{Completion of the proof of Theorem \ref{thm1}}\label{completion}

Lemma \ref{lem30} shows that $f$ has finite order $\rho(f)$. Thus it
remains only to prove that $f$ has finitely many poles and $f''$ has finitely many zeros, so assume that this is not the case.
Lemmas \ref{newlem0} and \ref{poleslem} give a positive real number
$d_1$ and $w \in \C$ with $|w| = r$ arbitrarily large, such that (\ref{lem10a}) holds on the disc
$B(w, d_1r)$. Let $\varepsilon$ and $K$ be  positive, with $\varepsilon$ small,
and let
$$
U_K = \left\{ z \in \C : \frac1K < |z| < K, \quad |z-1| > d_1 \right\} .
$$
Here $K$ is chosen so large that the harmonic measure with respect to $U_K$ satisfies 
\begin{equation}
 \label{hmest1}
\omega( z, S(0, 1/K) \cup S(0, K), U_K) < \varepsilon \quad \hbox{for} \quad
z \in U_K, \quad \frac12 < |z| < 2 .
\end{equation}
Denote by $d_j$ positive constants which are independent of $r$, $\varepsilon$, $K$ and $S$. 
Standard estimates from
\cite{Gun2} give a real number $S = S_r $ such that 
\begin{equation}
 \label{hmest2}
K < S < 2K \quad \hbox{and}  \quad | h(z) | \leq |z|^{d_2} \quad \hbox{for} \quad |z| = \frac{r}S \quad \hbox{and} \quad |z| = rS ,
\end{equation}
in which  $h(z) =   zf'''(z)/f''(z)$ as in Lemma \ref{wvlem} and $d_2 = \rho(f) + 1$. 
Let $w_1, \ldots, w_Q$ be the
poles of $h$ in $r/S \leq |z| \leq rS$. Applying Lemma \ref{lem2} with $L(r) = (\log r)^{1/2}$ shows that
$Q \leq d_3 (\log r)^{1/2} $.

On the annulus $A$ given by $r/S \leq |z| \leq rS$ set 
\begin{equation}
 \label{hmest3}
u(z) = \log |h(z)| - \log M_3  + \sum_{1 \leq j \leq Q} \log \frac{ |z-w_j| }{4Kr} \leq \log |h(z)| - \log M_3  , 
\end{equation}
where $M_3$ is as in (\ref{lem10a}) and may be assumed to be at least $1$, and the sum is empty if there are no poles $w_j$. 
Then $u$ is subharmonic  on $A$, with $u(z) \leq 0$ on the closure of $ B(w, d_1 r)$ by (\ref{lem10a}),  and 
\begin{equation}
 \label{hmest4}
u(z) \leq \log |h(z)| \leq d_2 \log |z| \leq d_2 \log (2Kr) \quad \hbox{for} \quad
z \in S(0, r/S) \cup S(0, rS) ,
\end{equation}
by  (\ref{hmest2}).
Hence (\ref{hmest1}) and the monotonicity of harmonic measure
yield
\begin{equation}
 \label{hmest5}
u(z) \leq \varepsilon  d_2 \log (2Kr)  \quad \hbox{for} \quad  \frac{r}2 < |z| < 2r.  
\end{equation}
Now Lemma \ref{cartan} shows that (\ref{cartanest}) holds, with $\Lambda = r/24$, 
for all $z$ outside a union $P_r $ of discs having sum of radii at most $r/4$.
Choose $s \in (r/2, 2r) \setminus E_0$, with $E_0$ as in Lemma \ref{wvlem}, such that
the circle $S(0, s)$ does not meet $P_r$. Thus Lemma \ref{wvlem} and (\ref{hmest5}) give rise to $\zeta_s \in S(0, s)$ such that
\begin{eqnarray*}
\frac13 \log s \leq \log |h(\zeta_s)|  &\leq&  \varepsilon  d_2 \log (2Kr) + \log M_3 +  \sum_{1 \leq j \leq Q} \log \frac{4Kr} { |\zeta_s-w_j| } \\
&\leq& \varepsilon  d_2 \log (2Kr) + \log M_3 +  Q \log (96K) \\
&\leq& \varepsilon (\rho(f)+1) \log (4Ks) + \log M_3  + d_3 (\log 2s)^{1/2} \log (96K) . 
\end{eqnarray*}
Since $\varepsilon$ may be chosen arbitrarily small, while $s$ is large,
this gives a contradiction and the proof of Theorem \ref{thm1} is complete. . 
\hfill$\Box$
\vspace{.1in}

\noindent
\textit{Remark.} Hypothesis (iii) on the multiplicities of poles may not be really essential 
for Theorem \ref{thm1} but it does play a key role in the above proof. If it is assumed merely that $f$ has finite lower order, 
then techniques such as P\'olya peaks should give annuli on which 
the analysis of Lemma \ref{lem2} can be applied, but it seems
difficult to ensure that these contain enough distinct poles of $f$ that the discs on which (\ref{lem10a}) holds are not so remote that
the method of Section \ref{completion} fails.

\noindent
School of Mathematical Sciences, University of Nottingham, NG7 2RD.\\
jkl@maths.nott.ac.uk

\end{document}